\tolerance=10000
\def\nagy{\font\caps=cmcsc10\caps}
\def\kis{\font\kisf=cmr5\kisf}

\voffset=.3 true in
\vsize=8.7 true in  

\newif\ifduplexpr
\duplexprfalse

\def\duplex{ \duplexprtrue \marginfalse}


\def\uj {\bigskip \rm}


\def\redef#1#2{\expandafter\ifx\csname #1\endcsname\relax
    \expandafter\edef\csname #1\endcsname{#2} 
    \else \message{redefinition of '\string#1'
    } \fi }

\newif\ifexist

\def\testfile#1{
\openin 0=#1
\ifeof 0
\existfalse
\else
\existtrue
\fi
\closein 0
}


\testfile{cite.inc}
\ifexist
\input cite.inc
\else
\message{!!!!!!!!!! One more pass needed for references !!!!!!!!!}
\fi


\immediate\openout 0=cite.inc


\immediate\openout 2=conten.inc
\def\writeitem#1#2{\write2
  {\string
  \line{#1\string\nagy{}#2\string\rm\string\leaderfill\string\quad\folio}}}


\def\boxit#1{\vbox{\hrule\hbox{\vrule#1\vrule}\hrule}}


\newif\ifmargin
\margintrue


\newcount\sorszam \sorszam=0
\def\sorszaminc {\advance\sorszam by1 
\ifnum \section=0 \else \the\section.\fi
\the\sorszam. }

\def\authorstr{}
\def\shorttitlestr{}

\def\author#1{\def\authorstr{#1} 
 \bigskip
 \centerline{#1}
}

\def\abstract#1{\bigskip
 \centerline{\nagy Abstract}\medskip
 \vbox{\narrower\narrower\noindent#1}
}

\def\shorttitle#1{\def\shorttitlestr{#1}}

\def\oldalszam{
\headline={
 \nagy
 \ifnum \folio > 1
 \ifodd \folio {\hfil \shorttitlestr \hfil\folio} 
 \else
 {\folio\hfil \authorstr \hfil}
 \fi
 \else \hfil
 \fi
}
}

\nopagenumbers
\oldalszam


\newcount\section \section=0
\def\newsection#1{
      \goodbreak
      \advance\section by1 \sorszam=0 
      \medskip\bigskip\centerline{\nagy \the\section. #1}
      \nobreak\nobreak\nobreak\nobreak
     \writeitem{\the\section.  }{#1}
     }
\newbox\labelboxx
\def\ugor{
\bigskip
\ifmargin
\hskip-1in {\box\labelboxx}
\vskip-\baselineskip
\fi
\noindent\bf}
\def\skippy{\enskip}

\def\theorem {    \ugor Theorem \sorszaminc \sl}   
\def\prop {       \ugor Proposition \sorszaminc \sl}
\def\lemma {      \ugor Lemma \sorszaminc \sl}      
\def\cor  {       \ugor Corollary \sorszaminc \sl} 
      
\def\proof { \medskip {\noindent \it Proof.\skippy}\rm}

\def\nullbox{\setbox0=\null \ht0=5pt \wd0=5pt \dp0=0pt \box0}
\def\eop {\hfill\hbox{\relax}\hfill \boxit{\nullbox}\goodbreak}    
\def\Eop {\vskip -\baselineskip \vskip -\belowdisplayshortskip \eop}


\def\label#1{  
\ifmargin
\vskip0in\hskip-1in {\kis#1}
\vskip-\baselineskip
\fi
\immediate
\write0{\string\redef{\string#1}{{\string\rm\the\section.\the\sorszam}}}}

\def\mark#1{  
\advance\sorszam by 1
\setbox\labelboxx=\hbox{\kis #1}
\immediate
\write0{\string\redef{\string#1}{{\string\rm
\ifnum \section=0 \else\the\section.\fi
\the\sorszam}}}
\advance\sorszam by -1
}

\def\labelref#1{ 
\immediate
\write0{\string\redef{\string#1}{{\string\rm\the\sorszam}}}}

\def\cite#1{\expandafter\ifx\csname#1\endcsname\relax
    ???
    \message{!!!!!!! Missing reference '\string#1' !!!!!!!}
    \else
\csname#1\endcsname
\fi}


\def\title#1{\shorttitle{#1}\centerline{\nagy #1}}


\def\references{
      \goodbreak
      \medskip\bigskip\centerline{\nagy References}
      \nobreak\nobreak
     \writeitem{\ }{References}
}

\def\i#1#2#3#4{\itemitem{\hbox to .5in{#1\hfil}}{#2, {\it #3}, #4}}




 \duplex

\title{Avoidable Sets in The Bicyclic Inverse Semigroup}

\author{N\'andor Sieben}

\medskip

\footnote{}{{2000 \it Mathematics subject classification:} 
20M18, 05C15} 

\footnote{}{{\it Keywords:\/}  avoidable set, additive partition, 
bicyclic inverse semigroup, bipartite graph} 

\abstract{A subset $U$ of a set $S$ with a binary operation 
is called {\it avoidable\/} if $S$ can be partitioned into two 
subsets $A$ and $B$ such that no element of $U$ can be 
written as a product of two distinct elements of $A$ or 
as the product of two distinct elements of $B$.  The 
avoidable sets of the bicyclic inverse semigroup are 
classified.  } 

\newsection{Introduction}

\uj If $(S,\cdot )$ is a set with a binary operation then a subset 
$U$ of $S$ is called {\it avoidable\/} if $S$ can be partitioned into 
two subsets $A$ and $B$ such that the partition {\it avoids} $U$, 
that is, no element of $U$ can be written as a product of 
two distinct elements of $A$ or as the product of two 
distinct elements of $B$.  Avoidable sets in $({\bf N},+)$ were 
first introduced by Alladi, Erd\H os and Hoggatt [AEH] and 
further studied in [Eva, Hog, HB, SZ, CL, ZC, Gra, Dum, 
De1].  Avoidable sets in groups were investigated in 
[De2].

In this paper we initiate the study of avoidable sets in 
inverse semigroups.  An inverse semigroup is a 
semigroup such that every element $s$ has a unique 
adjoint $s^{*}$ satisfying $ss^{*}s=s$ and $s^{*}ss^{*}=s^{*}$.  Inverse 
semigroups were first studied by Vagner [Vag] and 
Preston [Pr1] who considered inverse semigroups as the 
most promising class of semigroups for study.  An 
inverse semigroup is the next best thing to having an 
actual group.  While a group can be represented as 
bijections on a set, an inverse semigroup can be 
represented as partial bijections on a set.  In fact, a 
group is an inverse semigroup with a single idempotent.  
A comprehensive reference for inverse semigroups is 
[Pet].  

As a starting point of this study, we classify the 
maximal avoidable sets in the bicyclic inverse semigroup, 
which is perhaps the most important inverse semigroup.  
Its role in semigroup theory is similar to the role of ${\bf Z}$ 
in group theory.  It is one of the basic building blocks 
[Pr2] of the monogenic inverse semigroups, that is, 
inverse semigroups generated by single elements.  A 
possible continuation of our study could consider the 
other building blocks of monogenic inverse semigroups, in 
particular, the inverse semigroups generated by finite 
forward shifts.  

The {\it bicyclic inverse semigroup} ${\cal B}$ is the set 
$${\cal B}=\{(a,b)\mid a\ge 0,a+b\ge 0\}\subseteq {\bf Z}\times {\bf Z}$$
equipped with the following multiplication and inverse 
$$(a,b)(c,d)=(\max\{c+d,a\}-d,b+d),\qquad (a,b)^{*}=(a+b,-b).$$
Note that ${\cal B}$ can be represented as a semigroup of partial 
bijections of the nonnegative integers where the element 
$(a,b)$ is represented by the shift of the set $\{n\in {\bf Z}\mid 
n\ge a\}$ 
by $b$. Note that ${\cal B}$ is the inverse semigroup generated by 
the element $(0,1)$. 

Given $U\subseteq S$ the {\it associated graph} $G_{S,U}$ is the graph 
whose vertex set is $S$ and two vertices $r$ and $s$ are 
connected by an edge if $rs\in U$ or $sr\in U$.  Then $U$ is 
avoidable in $S$ exactly when $G_{S,U}$ is bipartite. So to show 
that a set $U$ is unavoidable in $S$, it suffices to find a 
cycle in $G_{S,U}$ with odd length.  

The author thanks the referee for suggestions that 
greatly simplified the paper.

\newsection{Even elements}

\uj An element $(a,b)$ of ${\cal B}$ is called {\it even\/} if $b$ is even, and 
{\it odd\/} if $b$ is odd.  Since the operation on ${\cal B}$ is written 
multiplicatively, it might be more appropriate, but less 
descriptive, to call even elements perfect squares.  It is 
clear that a product is odd exactly when the factors 
have different parity.  This shows that the set of odd 
elements is avoidable since the partition of ${\cal B}$ separating 
even and odd elements avoids the set off odd elements.  
Only special even elements of ${\cal B}$ can be in an avoidable 
set as the following proposition shows.  

\mark{evenprop}

\prop If $(a,b)\in {\cal B}$ is even such that $a\ge 1,$ $a+b\ge 
1$ and 
$b\neq 0$ then $U=\{(a,b)\}$ is unavoidable.  

\proof Let
$$p=\left(\alpha ,{b\over 2}\right),\quad q=\left(a,{b\over 2}\right
),\quad r=\left(a+{b\over 2},{b\over 2}\right)$$
where $\alpha$ will be chosen later.  Then $q\neq r$ since $b\neq 
0$.  
Also 
$$rq=\left(a+{b\over 2},{b\over 2}\right)\left(a,{b\over 2}\right
)=\left(\max\left\{a+{b\over 2},a+{b\over 2}\right\}-{b\over 2},b\right
)=(a,b)\in U.$$

If $b\ge 2$ then let $\alpha =a-1$. It is clear that $p\in {\cal B}$. Then
$$pq=\left(a-1,{b\over 2}\right)\left(a,{b\over 2}\right)=\left(\max\left
\{a+{b\over 2},a-1\right\}-{b\over 2},b\right)=(a,b)\in U,$$
$$rp=\left(a+{b\over 2},{b\over 2}\right)\left(a-1,{b\over 2}\right
)=\left(\max\left\{a-1+{b\over 2},a+{b\over 2}\right\}-{b\over 2}
,b\right)=(a,b)\in U.$$
It is clear that $p\neq q$, and since $b\neq -2$ we have $r\neq p$.  

If $b\le -2$ then let $\alpha =a+{b\over 2}-1$.  Then $a\ge -b\ge 
2$ and so 
$a+b\ge 0\ge -a+2$ which implies that $\alpha\ge 0$.  Since $a+b\ge 
1$, 
we also have $\alpha +{b\over 2}\ge 0$ and so $p\in {\cal B}$.  Then 
$$pq=\left(a+{b\over 2}-1,{b\over 2}\right)\left(a,{b\over 2}\right
)=\left(\max\left\{a+{b\over 2},a+{b\over 2}-1\right\}-{b\over 2}
,b\right)=(a,b),$$
and since $b\le -2$, we have $a+{b\over 2}\ge b-1$ which implies 
$$rp=\left(a+{b\over 2},{b\over 2}\right)\left(a+{b\over 2}-1,{b\over 
2}\right)=\left(\max\left\{a+b-1,a+{b\over 2}\right\}-{b\over 2},
b\right)=(a,b).$$
It is clear that $p\neq r$ and since $b\neq 2$ we have $p\neq q$.

In either case $\{p,q,r\}$ forms a triangle in $G_{{\cal B},U}$.
\eop

\uj The following figure shows the even elements of ${\cal B}$ 
that are not impossible in an avoidable set:
$$\matrix{&&&1_{{\cal B}}&\hbox{\rm {\it (0,2)}}&\hbox{\rm {\it (0,4)}}&
\ldots\cr
&&&\tt(1,0)&\times&\times&\cdots\cr
&&{\bf (}{\bf 2}{\bf ,}{\bf -}{\bf 2}{\bf )}&\tt(2,0)&\times&\times&\cr
&&\times&\tt(3,0)&\times&\times\cr
&{\bf (}{\bf 4}{\bf ,}{\bf -}{\bf 4}{\bf )}&\times&\tt(4,0)&\times&
\times\cr
&\times&\times&\tt(5,0)&\times&\times\cr
{\bf (}{\bf 6}{\bf ,}{\bf -}{\bf 6}{\bf )}&\times&\times&\tt(6,0)&
\times&\times&\cdots\cr
\times&\vdots&\vdots&\vdots&\vdots&\vdots\cr}
$$
This motivates the following notation:
$$\eqalign{{\cal D}&:=\{(a,-a)\mid a\ge 2,\hbox{\rm \ $a$ even}\}
,\quad {\cal E}:=\{(a,0)\mid a\ge 1\},\quad {\cal F}:=\{(0,b)\mid 
b\ge 2,\hbox{\rm \ $b$ even}\}.\cr}
$$
By Proposition~\cite{evenprop}, no even element outside 
${\cal D}\cup {\cal E}\cup {\cal F}\cup \{(0,0)\}$ can be in an avoidable set.  In the 
following sections we find the avoidable sets containing 
each type of these even elements.

\newsection{Sets containing $(a,-a)\in {\cal D}$}

\uj We investigate the possibility of $U$ being avoidable 
if $U\cap {\cal D}\neq\emptyset$.  

\mark{amalemma2}
\lemma If $(a,-a)\in {\cal D}$, $c\ge 0$ and $c+d\ge 1$ then 
$U=\{(a,-a),(c,d)\}\subseteq {\cal B}$ is unavoidable.  

\proof Let 
$$r=\left({a\over 2},-{a\over 2}\right),\quad s=\left({a\over 2}+
1,-{a\over 2}\right),\quad t=\left(c,d+{a\over 2}\right).$$
It is clear that $r,s\in {\cal B}$.  We also 
have $t\in {\cal B}$ because $c+d+a/2>c+d\ge 1>0$.  

It is clear that $r\neq s$.  We cannot have $r=t$ because 
that would imply $1\le c+d=a/2-a=-a/2<0$, which is 
impossible.  We cannot have $s=t$ either because that 
would imply $1\le c+d=a/2+1-a=1-a/2<1$. 

It is easy to see that $rs=(a,-a)\in U$. We 
also have
$$\eqalign{st&=\left({a\over 2}+1,-{a\over 2}\right)\left(c,d+{a\over 
2}\right)=\left(\max\left\{c+d+{a\over 2},{a\over 2}+1\right\}-d-{
a\over 2},d\right)=(c,d)\in U,\cr}
$$
and similar calculation shows that $rt=(c,d)\in U$. So 
$\{r,s,t\}$ forms a triangle in $G_{{\cal B},U}$. 
\eop

\uj We now consider the case when $c+d=0$, that is, 
$d=-c$.

\mark{amalemma4}

\lemma If $(a,-a)\in {\cal D}$, $c\neq a$ and ${a\over 2}<c$ then 
$U=\{(a,-a),(c,-c)\}\subseteq {\cal B}$ is unavoidable.  

\proof Since $a\neq c$ it is easy to see that
$$r=\left({a\over 2},-{a\over 2}\right),\quad s=\left({a\over 2}+
1,-{a\over 2}\right),\quad t=\left(c-{a\over 2},-c+{a\over 2}\right
)$$
are different elements of ${\cal B}$.

It is also easy to see that $rs=(a,-a)\in U$. We 
also have
$$ts=\left(c-{a\over 2},-c+{a\over 2}\right)\left({a\over 2}+1,-{
a\over 2}\right)=\left(\max\left\{1,c-{a\over 2}\right\}+{a\over 
2},-c\right)=(c,-c)\in U,$$
and similar calculation shows that $tr=(c,-c)\in U$ and so 
$\{r,s,t\}$ forms a triangle in $G_{{\cal B},U}$.
\eop

\mark{amacor3}
 
\cor If $(a,-a),(c,-c)\in {\cal D}$ and $a\neq c$ then 
$U=\{(a,-a),(c,-c)\}\subseteq {\cal B}$ is unavoidable.  

\proof Since $a\neq c$ we can assume, without loss of 
generality, that $a<c$. But then ${a\over 2}<c$ and so the result 
follows from the previous lemma. 
\eop

\mark{amaprop5}

\prop If $(a,-a)\in {\cal D}$, $0<c<e\le{a\over 2}$ and $c$, $e$ are odd then 
$U=\{(a,-a)$, $(c,-c)$, $(e,-e)\}$ is unavoidable.  

\proof It is clear that
$$r=\left(0,{{a-c-e}\over 2}\right),\quad s=\left({{a+c-e}\over 2}
,-{{a+c-e}\over 2}\right),\quad t=\left({{a-c+e}\over 2},-{{a-c+e}\over 
2}\right)$$
are different elements of ${\cal B}$. We have 
$$\eqalign{sr&=\left({{a+c-e}\over 2},-{{a+c-e}\over 2}\right)\left
(0,{{a-c-e}\over 2}\right)\cr
&=\left(\max\left\{{{a-c-e}\over 2},{{a+c-e}\over 2}\right\}-{{a-
c-e}\over 2},-c\right)=(c,-c)\in U,\cr}
$$
$$\eqalign{st&=\left({{a+c-e}\over 2},-{{a+c-e}\over 2}\right)\left
({{a-c+e}\over 2},-{{a-c+e}\over 2}\right)\cr
&=\left(\max\left\{0,{{a+c-e}\over 2}\right\}+{{a-c+e}\over 2},-a\right
)=(a,-a)\in U,\cr}
$$
$$\eqalign{tr&=\left({{a-c+e}\over 2},-{{a-c+e}\over 2}\right)\left
(0,{{a-c-e}\over 2}\right)\cr
&=\left(\max\left\{{{a-c-e}\over 2},{{a-c+e}\over 2}\right\}-{{a-
c-e}\over 2},-e\right)=(e,-e)\in U,\cr}
$$
and so $\{r,s,t\}$ forms a triangle in $G_{{\cal B},U}$.
\eop

\mark{amaprop6}

\prop If $(a,-a)\in {\cal D}$ and $0<c<{a\over 2}$ then $U=\{(a,-
a),(c,-c),(0,0)\}$ is 
unavoidable.

\proof It is clear that
$$p=(0,0),\quad q=(c,-c),\quad r=(a-c,-(a-c)),\quad s=(0,a-c),\quad 
t=(a,-a)$$
are different elements of ${\cal B}$.  It is easy to check that 
$pq=(c,-c)$, $qr=(a,-a)$, $rs=(0,0)$, $ts=(c,-c)$ and 
$pt=(a,-a)$. Thus $\{p,q,r,s,t\}$ forms a cycle with odd 
length in $G_{{\cal B},U}$.  \eop

\uj We are going to denote the remainder of $y$ modulo $m$ 
by $[y]_m$.

\mark{amacol}

\prop Let $(a,-a)\in {\cal D}$ and $c$ be odd.  If $0<c<{a\over 2}$ then 
$U=\{(a,-a),(c,-c)\}$ is avoidable.  If $c={a\over 2}$ then 
$U=\{(a,-a),\left({a\over 2},{a\over 2}\right),(0,0)\}$ is avoidable.  

\proof Let 
$$A=\left\{(x,y)\mid{{a-2c}\over 2}\le [y]_{a-c}\le{{2a-3c-1}\over 
2}\right\}\setminus\left\{\left({a\over 2}+k(a-c),-{a\over 2}-k(a
-c)\right)\mid k=0,1,\ldots\right\}$$
and $B={\cal B}\setminus A$.  Note that $[v]_{a-c}={{a-2c}\over 2}$ if and only if 
$v=-{a\over 2}+l(a-c)$ for some $l\in {\bf Z}$.  We show that the 
partition $\{A,B\}$ avoids $U$.  The following figure shows the 
partition when $a=8$ and $c=3$.  Note that 
${{a-2c}\over 2}=\tt1$ and ${{2a-3c-1}\over 2}=\tt3$ in this case.  
$$\matrix{&&-9&-8&-7&-6&-5&-{a\over 2}&-3&-2&-1&0&\tt1&2&\tt3&4&a
-c&\cdots\cr
\cr
0&&.&.&.&.&.&.&.&.&.&1&0&0&0&1&1\cr
1&&.&.&.&.&.&.&.&.&1&1&0&0&0&1&1\cr
2&&.&.&.&.&.&.&.&0&1&1&0&0&0&1&1\cr
3&&.&.&.&.&.&.&0&0&1&1&0&0&0&1&1\cr
-a/2&&.&.&.&.&.&{\bf 1}&0&0&1&1&0&0&0&1&1\cr
5&&.&.&.&.&1&0&0&0&1&1&0&0&0&1&1\cr
6&&.&.&.&1&1&0&0&0&1&1&0&0&0&1&1\cr
7&&.&.&0&1&1&0&0&0&1&1&0&0&0&1&1\cr
-a&&.&0&0&1&1&0&0&0&1&1&0&0&0&1&1\cr
9&&{\bf 1}&0&0&1&1&0&0&0&1&1&0&0&0&1&1\cr
10&&0&0&0&1&1&0&0&0&1&1&0&0&0&1&1\cr
\vdots\cr}
$$
If $(\alpha ,-\alpha )=(x,y)(w,z)=(\max\{w+z,x\}-z,y+z)$ then either 
$w+z\ge x$ and $w=\alpha$ or $x\ge w+z$ and $x-z=\alpha$.  In the 
first case we have $0\le x+y\le w+z+y=\alpha -\alpha =0$ which 
can only happen if $y=-x$.  In the second case we have 
$y=-\alpha -z=-x$ as well.  Thus we only have to show that 
$(x,-x)(w,z)$ is not in $U$ unless the two factors are in 
separate classes.  

Let $s=(x,y)\neq t=(w,z)$ and suppose that $s$ and $t$ are in 
the same equivalence class.  If $st\in \{(a,-a),(c,-c)\}$ then 
either $y+z=-a$ or $y+z=-c$.  So we need to study the 
effect of the maps $v\mapsto -v-a$ and $v\mapsto -v-c$ on the 
congruence classes modulo $a-c$.  Since $a$ and $c$ are 
congruent modulo $a-c$, we only need to study one of the 
maps.  If ${{a-2c}\over 2}\le [v]_{a-c}\le a-2c$ then 
$$0=-(a-2c)-a+2(a-c)\le [-v-a]_{a-c}\le -{{a-2c}\over 2}-a+2(a-c)
={{a-2c}\over 2}.$$
If $a-2c<[v]_{a-c}\le{{2a-3c-1}\over 2}$ then 
$${{2a-3c+1}\over 2}=-{{2a-3c-1}\over 2}-a+3(a-c)\le [-v-a]_{a-c}
<-(a-2c)-a+3(a-c)=a-c.$$
So 
$$[v]_{a-c}\in\left[{{a-2c}\over 2},{{2a-3c-1}\over 2}\right]\Leftrightarrow 
[-v-a]_{a-c},[-v-c]_{a-c}\in\left[0,{{a-2c}\over 2}\right]\cup\left
[{{2a-3c+1}\over 2},a-c\right).$$

First consider the case when $c<{a\over 2}$.  We must have 
$[y]_{a-c}={{a-2c}\over 2}=[z]_{a-c}$ and so $y=-{a\over 2}+k(a-c
)$ and 
$z=-{a\over 2}+l(a-c)$ for some $k,l\in {\bf Z}$.  Since 
$0\le x=-y={a\over 2}-k(a-c)$, we must have $k\le 0$ and so $s\in 
B$.  

If $st=(c,-c)$ then $-c=y+z=-a+(k+l)(a-c)$ and so 
$k+l=1$.  Hence $l=1-k\ge 1$ and so $t\in A$ which is a 
contradiction.  If $st=(a,-a)$ then 
$-a=y+z=-a+(k+l)(a-c)$ and so $k+l=0$.  If $k=0$ 
then $l=0$ and so $z=y$.  Since $s\neq t$ we must have 
$w\neq x=-y=-z$ and so $t\in A$ which is a contradiction.  If 
$k<0$ then $l>0$ and so $t\in A$ which is again a 
contradiction.  

Next consider the case when $c={a\over 2}$.  If $st\in \{(c,-c),(
a,a)\}$ 
then again we must have $[y]_{a-c}={{a-2c}\over 2}=[z]_{a-c}$ and so 
$y=k{a\over 2}$ and $z=l{a\over 2}$ for some $k,l\in {\bf Z}$.  Since 
$0\le x=-y=-k{a\over 2}$, we have $k\le 0$.  If $k<0$ then $s\in 
B$ and 
we get a contradiction like we did in the $c<{a\over 2}$ case.  If 
$k=0$ then $s=(0,0)\in A$ and so $t=(c,-c)\in B$ or 
$t=(a,-a)\in B$ which is a contradiction.  

If $st=(0,0)$ then we must have $s=(x,-x)$ and $t=(0,x)$ 
and since $s\neq t$, we also know that $x\ge 1$.  Note that in 
this case $a-c=c,$ ${{a-2c}\over 2}=0$ and ${{2a-3c-1}\over 2}={{
a-2}\over 4}$.  If $a=2$ 
then $B=\{(x,-x)\mid x=1,2,\ldots \}$ and so $s\in B$ while $t\in 
A$.  
Since $c$ is odd we cannot have $a=4$.  If $a\ge 6$ then 
$0<{{a-2}\over 4}<{a\over 2}=a-c$.  So we have
$$0\le [v]_{a-c}\le{{a-2}\over 4}$$
if and only if either 
$$[-v]_{a-c}\ge -{{a-2}\over 4}+(a-c)={{a+2}\over 4}$$
or $[-v]_{a-c}=0$.  So since $s$ and $t$ are in the same 
equivalence class we must have $-x=-{a\over 2}k$ for some 
positive $k$ and so $s\in B$ while $t\in A$.  \eop

\mark{amaonly}

\prop If $U\cap {\cal D}\neq\emptyset$ and $U$ is maximal avoidable then $
U$ is 
one of the avoidable sets of Proposition~\cite{amacol}.  

\proof Let $(a,-a)\in U\cap {\cal D}$.  If $(a,-a)\neq (x,y)\in U$ then by 
Lemma~\cite{amalemma2}, we must have $y=-x$.  By 
Lemma~\cite{amalemma4} and Corollary \cite{amacor3} we 
know that $x$ cannot be even and $x\le{a\over 2}$.  If $0<y<{a\over 
2}$ then 
by Propositions~\cite{amaprop5}, \cite{amaprop6} and 
\cite{amacol} $U=\{(a,-a),(x,-x)\}$ is maximal avoidable.  If 
$y=0$ or $y={a\over 2}$ then by Propositions~\cite{amaprop5} and 
\cite{amaprop6}, $U$ cannot have yet another element 
$(w,-w)$ unless $w=0$ or $w={a\over 2}$.  This fact and 
Proposition~\cite{amacol} implies that 
$U=\{(a,-a),({a\over 2},-{a\over 2}),(0,0)\}$ is maximal avoidable.  We 
considered all the possibilities so these are the only 
maximal avoidable sets intersecting ${\cal D}$.  \eop

\newsection{Sets containing $(0,b)\in {\cal F}$}

\uj We investigate the possibility of $U$ being avoidable if 
$U\cap {\cal F}\neq\emptyset$.  Our main tool is the fact that ${\cal F}
={\cal D}^{*}$, which 
allows us to transform the results of the previous 
section.  

\prop The set $U$ is avoidable if and only if $U^{*}$ is 
avoidable. Furthermore, $U$ is maximal avoidable if and 
only if $U^{*}$ is maximal avoidable.

\proof First, assume $U$ is avoidable.  Then $U$ can be 
partitioned into two subsets $A$ and $B$ such that the 
partition avoids $U$.  Now $\{A^{*},B^{*}\}$ is a partition of $U^{
*}$.  
If $x,y\in A^{*}$ and $x\neq y$ then $x^{*},y^{*}\in A$ and so $y^{
*}x^{*}\notin U$, 
which means $xy=(y^{*}x^{*})^{*}\notin U^{*}$.  Similar argument shows 
that no element of $U^{*}$ is the product of two different 
elements of $B^{*}$.  

Now if $U^{*}$ is avoidable then $U=U^{**}$ is also avoidable by 
the previous argument.  

The second part of the proposition follows from the fact 
that if $U$ and $V$ are subsets of ${\cal B}$ then $U\subseteq V$ exactly 
when $U^{*}\subseteq V^{*}$.
\eop

\mark{0bprop9}

\prop If $U\cap {\cal F}\neq\emptyset$ and $U$ is maximal avoidable then either 
$U=\{(0,b),(0,d)\}$ where $(0,b)\in {\cal F}$, $d$ is odd and $0<
d<{b\over 2}$, or
$U=\{(0,b),\left(0,{b\over 2}\right),(0,0)\}$ where $(0,b)\in {\cal F}$ and ${
b\over 2}$ is odd.

\proof By the previous proposition, $U$ is maximal 
avoidable exactly when $U^{*}$ is maximal avoidable. Since 
$U\cap {\cal F}\neq\emptyset$ and ${\cal F}^{*}={\cal D}$, we must have $
U^{*}\cap {\cal D}\neq\emptyset .$ Hence $U^{*}$ 
is one of the maximal avoidable sets of 
Proposition~\cite{amacol}. 
\eop

\newsection{Sets containing $(a,0)\in {\cal E}$}

\uj We investigate the possibility of $U$ being avoidable if 
$U\cap {\cal E}\neq\emptyset$.  

\mark{a0blemma}

\lemma If $(a,0)\in {\cal E}$ and $\max\{c,c+d\}\ge a$ then 
$U=\{(a,0),(c,d)\}\subseteq {\cal B}$ is unavoidable.  

\proof If $c\ge a$ then $(c,d)(a,0)=(c,d)$.  If $c+d\ge a$ then 
$(a,0)(c,d)=(c,d)$.  In either case $\{\left(0,0\right),\left(a,0\right
),\left(c,d\right)\}$ 
forms a triangle in $G_{{\cal B},U}$.   
\eop

\mark{a0cor}

\cor If $a\neq c$ then $U=\{(a,0),(c,0)\}\subseteq {\cal E}$ is unavoidable.

\proof Without loss of generality we can assume that 
$a<c$ and so the result follows from 
Lemma~\cite{a0blemma}.  
\eop

\mark{a0lem}

\lemma If $(a,0)\in {\cal E}$, $d$ and $f$ are odd, $d\neq f$, 
$\max\{c,c+d,e,e+f\}<a$ then $U=\{(a,0),(c,d),(e,f)\}\subseteq {\cal B}$  is 
unavoidable.  

\proof Without loss of generality we can assume that 
$d<f$.  Let $x={{f-d}\over 2}$, $y={{f+d}\over 2}$ and 
$$p=(\alpha ,x),\quad q=(\beta ,y),\quad r=(c+y,-x),\quad s=(a,x)
,\quad t=(a,-x).$$
where $\alpha$ and $\beta$ will be chosen later.  Since $y-x=d$ and 
$x+y=f$ are odd and so not zero, we have $p\neq q\neq r$.  
Since $x>0$, we also have $r\neq s\neq t\neq p$.  

Since 
$$c+y=c+{{f+d}\over 2}={{c+f+c+d}\over 2}\ge{{c+d+c+d}\over 2}=c+
d\ge 0$$
and $c+y-x=c+d\ge 0$, we have $r\in {\cal B}$.  Since 
$a+x>a\ge 0$, we have $s\in {\cal B}$.  We have $t\in {\cal B}$ because 
$$a-x=a-{{f-d}\over 2}={{(a-f)+a+d}\over 2}>{{e+a-c}\over 2}>{e\over 
2}\ge 0.$$
Since $a+x>c+d+x=c+y$ we have
$$rs=(c+y,-x)(a,x)=(\max\{a+x,c+y\}-x,0)=(a,0)\in U.$$
Also 
$$ts=(a,-x)(a,x)=(\max\{a+x,a\}-x,0)=(a,0)\in U.$$

First, assume that $e+x<a$.  Let $\alpha =e$ and 
$\beta =\min\{c,e+x\}\ge 0$.  Then clearly $p\in {\cal B}$ and since 
$$\eqalign{\min\{c,e+x\}+y&=\min\{c+y,e+x+y\}=\min\left\{{{c+f+c+
d}\over 2},e+f\right\}\cr
&\ge\min\{c+d,e+f\}\ge 0,\cr}
$$
we must have $q\in {\cal B}$. Now we have
$$pt=(e,x)(a,-x)=(\max\{a-x,e\}+x,0)=(a,0)\in U,$$
$$qp=(\min\{c,e+x\},y)(e,x)=(\max\{e+x,\min\{c,e+x\}\}-x,y+x)=(e,
f)\in U,$$
$$rq=(c+y,-x)(\min\{c,e+x\},y)=(\max\{\min\{c,e+x\}+y,c+y\}-y,d)=
(c,d)\in U,$$
and so $\{p,q,r,s,t\}$ forms a cycle with odd length in 
$G_{{\cal B},U}$.  

Next, assume that $e+x\ge a$.  We need to consider two 
subcases.  In the first subcase we assume that $e<c$.  
Let $\alpha =0$ and $\beta =e$.  It is clear that $p\in {\cal B}$.  Since 
$e+y=e+x+y-x\ge a+d>c+d\ge 0$, we have $q\in {\cal B}$. It is 
clear that $pt=(0,x)(a,-x)=(a,0)\in U$ and 
$pq=(0,x)(e,y)=(e,f)\in U$. We also have
$$rq=(c+y,-x)(e,y)=(\max\{e+y,c+y\}-y,d)=(c,d)\in U,$$
and so $\{p,q,r,s,t\}$ forms a cycle with odd length in 
$G_{{\cal B},U}$.  

In the second subcase we assume that $e\ge c$.  This 
implies that $e+d\ge c+d\ge 0$.  Let $\alpha =e+y$ and $\beta =c$.  
Since 
$$e+y=e+{{f+d}\over 2}={{e+f+e+d}\over 2}\ge 0$$
and $e+y+x=e+f\ge 0$, we have $p\in {\cal B}$. Since 
$$c+y=c+{{f+d}\over 2}={{c+f+c+d}\over 2}>{{c+d+c+d}\over 2}=c+d\ge 
0,$$
we have $q\in {\cal B}$. Since $a-x>e+f-x=e+y+x-x=e+y$ 
we have
$$pt=(e+y,x)(a,-x)=(\max\{a-x,e+y\}+x,0)=(a,0)\in U.$$
Also  
$$pq=(e+y,x)(c,y)=(\max\{c+y,e+y\}-y,f)=(e,f)\in U,$$
$$rq=(c+y,-x)(c,y)=(\max\{c+y,c+y\}-y,d)=(c,d)\in U,$$
and so $\{p,q,r,s,t\}$ forms a cycle with odd length in 
$G_{{\cal B},U}$.  
\eop

\mark{a0prop}

\prop If $(a,0)\in {\cal E}$, $d$ is odd and $d<a$ then 
$$U=\{(a,0),(0,0)\}\cup \{(c,d)\in {\cal B}\mid\max\{c,c+d\}<a\}$$
is avoidable.

\proof First assume that $d>0$. For $y\in {\bf Z}\setminus \{0\}$ define 
$$\eqalign{\phi (y)&=\cases{0&if $0<[y]_d\le{{d-1}\over 2}$\cr
0&if $[y]_d=0$ and $y<0$\cr
1&otherwise.\cr}
\cr}
$$
Note that if $y\neq 0$ then $\phi (y)\neq\phi (-y)$.  Let 
$$A=\{(x,y)\in {\cal B}\mid y\neq 0\hbox{\rm \ and }\phi (y)=0\}\cup 
\{(x,y)\in {\cal B}\mid y=0\hbox{\rm \ and }x<a\}$$
and $B={\cal B}\setminus A$.  We show that the partition $\{A,B\}$ 
avoids $U$.  The following figure shows the partition 
when $a=6$ and $d=5$:  
$$\matrix{&&-7&-6&-d&-4&-3&-2&-1&0&1&{{d-1}\over 2}&3&4&d&6&7&\cdots\cr
\cr
0&&.&.&.&.&.&.&.&\tt0&0&0&1&1&{\bf 1}&0&0\cr
1&&.&.&.&.&.&.&1&\tt0&0&0&1&1&{\bf 1}&0&0\cr
2&&.&.&.&.&.&1&1&\tt0&0&0&1&1&{\bf 1}&0&0\cr
3&&.&.&.&.&0&1&1&\tt0&0&0&1&1&{\bf 1}&0&0\cr
4&&.&.&.&0&0&1&1&\tt0&0&0&1&1&{\bf 1}&0&0\cr
5&&.&.&{\bf 0}&0&0&1&1&\tt0&0&0&1&1&{\bf 1}&0&0\cr
a&&.&1&{\bf 0}&0&0&1&1&\tt1&0&0&1&1&{\bf 1}&0&0\cr
7&&1&1&{\bf 0}&0&0&1&1&\tt1&0&0&1&1&{\bf 1}&0&0\cr
\vdots\cr}
$$
Let $s=(x,y)\neq t=(w,z)$ and suppose that $s$ and $t$ are in 
the same equivalence class.  If $st=(u,0)$ then $y+z=0$ 
and so $y=0=z$.  Hence if $st=(0,0)$ then we must have 
$s=(0,0)=t$ which is a contradiction.  If $st=(a,0)$ then 
without loss of generality we can assume that $x>w$ 
and so $x=a$ and $w<a$.  Thus $s\in B$ and $t\in A$ which is a 
contradiction.  

Now assume $(c,d)\in U$ and $(c,f)$ can be written as $st$ or 
$ts$.  Then $y+z=d$ and so we either have $y=0$ and 
$z=d$ or we have $y=d$ and $z=0$.  If $y=0$ and $z=d$ 
then $t\in B$ and so we must have $s\in B$ which implies that 
$x\ge a$.  Hence $(c,d)=st=(x,0)(w,d)=(\max\{w+d,x\}-d,d)$ 
and so $c=\max\{w+d,x\}-d\ge x-d\ge a-d$ which is a 
contradiction.  If $y=d$ and $z=0$ then $s\in B$ and so we 
must have $t\in B$ which implies $w\ge a$.  Hence 
$(c,d)=(x,d)(w,0)=(\max\{w,x\},d)$ and so 
$c=\max\{w,x\}\ge w\ge a$ which is a contradiction.  

In case $d<0$ the proof is similar but we need to replace
the definition of $\phi$ by the following:
$$\eqalign{\phi (y)&=\cases{0&if $0<[y]_d\le{{d-1}\over 2}$\cr
0&if $[y]_d=0$ and $y>0$\cr
1&otherwise.\cr}
\cr}
$$
\Eop

\mark{a0prop5}

\prop If $U\cap {\cal E}\neq\emptyset$ and $U$ is maximal avoidable then $
U$ is 
the avoidable set of Proposition~\cite{a0prop}.  

\proof Let $(a,0)\in U\cap {\cal E}$.  If $(a,0)\neq (c,d)\in U$ then by 
Lemma~\cite{amalemma2}, $(c,d)\notin {\cal D}$, by 
Lemma~\cite{0bprop9}, $(c,d)\notin {\cal F}$ and by 
Corollary~\cite{a0cor}, $(c,d)\notin {\cal E}$.  If $(c,d)\neq (0
,0)$ then $d$ is 
odd and by Lemma~\cite{a0blemma}, $\max\{c,c+d\}<a$.  If 
$(w,z)$ is also in $U$ and $(c,d)\neq (w,z)\neq (0,0)$ then again 
$\max\{w,w+z\}<a$ and by Lemma~\cite{a0lem}, $z=d$.  \eop

\newsection{Sets containing $(0,0)$}

\mark{00prop1}

\prop If $d$ and $f$ are odd and $d\neq f$ then 
$U=\{(0,0),(c,d),(e,f)\}\subseteq {\cal B}$ is unavoidable. 

\proof Without loss of generality we can assume that 
$d<f$.  Let $x={{f-d}\over 2}$ and $y={{f+d}\over 2}$.  Note that $
x\neq y\neq -x$ 
since $y-x=d$ and $x+y=f$ are odd.  Also note that 
$c+y=c+d+x\ge x\ge 0$. 

First, we consider the case when $c\le e$.  If 
$$p=(0,x),\quad q=(x,-x),\quad r=(c,y),\quad s=(e+y,x),\quad t=(e
,y)$$
then $p\neq q\neq r\neq s\neq t\neq p$.  Since $c+y\ge 0$, $r\in 
{\cal B}$.  Since 
$$e+y={{e+f+e+d}\over 2}\ge{{e+f+c+d}\over 2}\ge 0$$
we have $s,t\in {\cal B}$.  Since $e+y+x=e+f\ge 0$ we have $s\in 
{\cal B}$.  
It is easy to see that $qp=(0,0)\in U$ and $pt=(e,f)\in U$.  
We have 
$$qr=(x,-x)(c,y)=(\max\{c+y,x\}-y,y-x)=(c,d)\in U,$$
$$sr=(e+y,x)(c,y)=(\max\{c+y,e+y\}-y,f)=(e,f)\in U,$$
$$st=(e+y,x)(e,y)=(\max\{e+y,e+y\}-y,f)=(e,f)\in U$$
and so $\{p,q,r,s,t\}$ forms a cycle with odd length in $G_{{\cal B}
,U}$.

Next, we consider the case when $e<c$.  We need to 
consider two subcases.  In the first subcase we assume 
that $e+y\ge 0$.  It is clear that 
$$p=(0,x),\quad q=(x,-x),\quad r=(c,y),\quad s=(c+y,-x),\quad t=(
e,y)$$
are in ${\cal B}$ and $p\neq q\neq r\neq s\neq t\neq p$.  As before, we have 
$qp=(0,0)\in U$, $pt=(e,f)\in U$ and $rq=(c,d)\in U$.  We also 
have 
$$qr=(x,-x)(c,y)=(\max\{c+y,x\}-y,d)=(c,d)\in U,$$
$$sr=(c+y,-x)(c,y)=(\max\{c+y,c+y\}-y,d)=(c,d)\in U,$$
and so $\{p,q,r,s,t\}$ forms a cycle with odd length in $G_{{\cal B}
,U}$.  

In the second subcase we assume that $e+y<0$. If
$$p=(0,x),\quad q=(e+x,y),\quad r=(c+y,-x),\quad s=(c-x,y),\quad 
t=(x,-x)$$
then $p\neq q\neq r\neq s\neq t\neq p$.  It is clear that $p,q,t\in 
{\cal B}$.  
Since $c+y\ge x\ge 0$ and 
$c-x=c+y-x-y=c+d-y>c+d+e\ge e\ge 0$, we have 
$r,s\in {\cal B}$.  It is easy to see that $tp=(0,0)\in U$ and 
$qp=(e,f)\in U$.  Since $c\ge e+x$ we have 
$$rq=(c+y,-x)(e+x,y)=(\max\{e+x+y,c+y\}-y,d)=(c,d)\in U.$$
We also have 
$$rs=(c+y,-x)(c-x,y)=(\max\{c-x+y,c+y\}-y,d)=(c,d)\in U,$$
$$st=(c-x,y)(x,-x)=(\max\{0,c-x\}+x,d)=(c,d)\in U,$$
and so $\{p,q,r,s,t\}$ forms a cycle with odd length in $G_{{\cal B}
,U}$.  
\eop

\mark{00prop2}

\prop If $d$ is odd then $U=\{(0,0)\}\cup \{(e,f)\in {\cal B}\mid 
f=d\}$ is 
avoidable.  

\proof Let 
$$A=\left\{(x,y)\mid [x]_{|d|}\le{{|d|-1}\over 2}\right\}\setminus 
\{y=kd\mid k=1,2,\ldots \}$$
and $B={\cal B}\setminus A$. We show that the partition $\{A,B\}$ avoids $
U$. 
The following figure shows the partition when $d=5$. 
Note that ${{|d|-1}\over 2}=\tt2$ in this case.
$$\matrix{&&-7&-6&-5&-4&-3&-2&-1&0&1&\tt2&3&4&d&6&7&8&9&2d&11&12&
13&\cdots\cr
\cr
0&&.&.&.&.&.&.&.&0&0&0&1&1&{\bf 1}&0&0&1&1&{\bf 1}&0&0&1\cr
1&&.&.&.&.&.&.&1&0&0&0&1&1&{\bf 1}&0&0&1&1&{\bf 1}&0&0&1\cr
2&&.&.&.&.&.&1&1&0&0&0&1&1&{\bf 1}&0&0&1&1&{\bf 1}&0&0&1\cr
3&&.&.&.&.&0&1&1&0&0&0&1&1&{\bf 1}&0&0&1&1&{\bf 1}&0&0&1\cr
4&&.&.&.&0&0&1&1&0&0&0&1&1&{\bf 1}&0&0&1&1&{\bf 1}&0&0&1\cr
5&&.&.&0&0&0&1&1&0&0&0&1&1&{\bf 1}&0&0&1&1&{\bf 1}&0&0&1\cr
6&&.&1&0&0&0&1&1&0&0&0&1&1&{\bf 1}&0&0&1&1&{\bf 1}&0&0&1\cr
7&&1&1&0&0&0&1&1&0&0&0&1&1&{\bf 1}&0&0&1&1&{\bf 1}&0&0&1\cr
\vdots\cr}
$$
Let $s=(x,y)\neq t=(w,z)$ and suppose that $s$ and $t$ are in 
the same equivalence class.  If $st=(0,0)$ then $s=(x,-x)$ 
and $t=(0,x)$.  But $s$ and $t$ cannot be in the same 
equivalence class unless $y=0=z$.  So $s=t$ which is a 
contradiction.  If $st\neq (0,0)$ then either we have $y=0$ 
and $z=d$ or we have $y=d$ and $z=0$.  Since $y=-x$ and 
so negative, $y=0$ and $z=d$ and so $s\in A$ but $t\in B$ which 
is a contradiction.  \eop

\mark{00prop3}

\prop If $d$ is odd then $U=\{(0,0)\}\cup \{(e,f)\in {\cal B}\mid 
f=d\}$ is 
maximal avoidable.  

\proof Clearly $U$ has infinitely many elements.  In 
particular, if $d>0$ then $(0,d)\in U$ and if $d<0$ then 
$(-d,d)\in U$.  If $(x,y)\in U$ and $y$ is odd then by 
Proposition~\cite{00prop1} $y=d$.  If $(x,y)\neq (0,0)$ but $y$ is 
even then $(x,y)\in {\cal D}\cup {\cal F}\cup {\cal E}$.  This is impossible because 
then Propositions~\cite{amaonly}, \cite{0bprop9} and 
\cite{a0prop} would imply that $U$ is finite.  \eop

\newsection{Conclusion}

\uj We are in position to give a full classification of the 
maximal avoidable sets of ${\cal B}$.

\theorem The maximal avoidable sets in ${\cal B}$ are the 
following:  
\item{{\rm(a)}}{$\{(a,b)\in {\cal B}\mid b\hbox{\rm \ is odd}\}$; } 
\item{{\rm(b)}}{$\{(0,0)\}\cup \{(e,f)\in {\cal B}\mid f=d\}$ where $
d$ is a 
fixed odd number; } 
\item{{\rm(c)}}{$\{(a,0),(0,0)\}\cup \{(c,d)\in {\cal B}\mid\max\{
c,c+d\}<a\}$ 
where $a$ and $d$ are fixed such that $d<a>0$ and $d$ is odd; } 
\item{{\rm(d)}}{$\{(a,-a),(c,-c)\}$ where $a>0$, $a$ is even, 
$0<c<{a\over 2}$ and $c$ is odd; } 
\item{{\rm(e)}}{$\{(a,-a),\left({a\over 2},-{a\over 2}\right),(0,
0)\}$ where $a>0$, $a$ is 
even and ${a\over 2}$ is odd; } 
\item{{\rm(f)}}{$\{(0,b),(0,d)\}$ where $b>0$, $b$ is even, 
$0<d<{b\over 2}$ and $d$ is odd; } 
\item{{\rm(g)}}{$\{(0,b),\left(0,{b\over 2}\right),(0,0)\}$ where $
b>0$, $b$ is even 
and ${b\over 2}$ is odd.  } 

\proof Suppose $U$ is maximal avoidable.  If $U$ has no even 
elements then $U$ is the set of part (a).  So suppose $U$ 
has an even element $s$.  Then by 
Proposition~\cite{evenprop}, $s=(0,0)$ or $s\in {\cal E}\cup {\cal D}
\cup {\cal F}$.  
If $s\in {\cal E}\cup {\cal D}\cup {\cal F}$ then by Propositions~\cite{amaonly}, 
\cite{0bprop9} and \cite{a0prop}, $U$ is one of the sets in 
parts (c,d,e,f,g). Assume $s=(0,0)$ and $s\neq t\in U$. If $t$ is 
even then again $t\in {\cal E}\cup {\cal D}\cup {\cal F}$ and so $
U$ is one of the sets 
in parts (c,e,g). If $U$ does not have any more even 
elements then by Proposition~\cite{00prop1}, $U$ is the set 
in part (b). 
\eop

\uj The notion of avoidable sets can be generalized.  We 
could call a set $U$ of $S$ $n$-avoidable if there is a 
partition of $S$ into $n$ subsets such that no element of $U$ 
can be written as a product of two distinct elements of 
the same subset.  It would be interesting to know if 
there are sets that are not 3-avoidable in a 
group or an inverse semigroup.

\references

\uj
\i{[AEH]}{K. Alladi, P. Erd\H os, V.E. Hoggatt}
 {On additive partitions of integers}
 {Discrete Math. {\bf 22} (1978), no. 3, 201--211.}

\i{[CL]}{T.Y. Chow, C.D. Long}
 {Additive partitions and continued fractions}
 {Ramanujan J. {\bf 3} (1999), no. 1, 55--72.}

\i{[De1]}{M. Develin}
 {A complete categorization of when generalized Tribonacci sequences can be avoided by additive partitions}
 {Electron. J. Combin. {\bf7} (2000), no. 1, Research Paper 53, 7 pp.}

\i{[De2]}{M. Develin}
 {Avoidable sets in groups}
 {Ars Combin. {\bf 65} (2002), 279--297.}

\i{[Dum]}{I. Dumitriu}
 {On generalized Tribonacci sequences and additive partitions}
 {Discrete Math. {\bf 219} (2000), no. 1-3, 65--83.}

\i{[Eva]}{R.J. Evans}
 {On additive partitions of sets of positive integers}
 {Discrete Math. {\bf 36} (1981), no. 3, 239--245.}

\i{[Gra]}{D. Grabiner}
 {Continued fractions and unique additive partitions}
 {Ramanujan J. {\bf 3} (1999), no. 1, 73--81.}

\i{[Hog]}{V. Hoggatt}
 {Additive partitions of the positive integers}
 {Fibonacci Quart. {\bf 18} (1980), no. 3, 220--226.}

\i{[HB]}{V. Hoggatt, M. Bicknell-Johnson}
 {Additive partitions of the positive integers and generalized 
  Fibonacci representations}
 {Fibonacci Quart. {\bf 22} (1984), no. 1, 2--21.}

\i{[Pet]}{M. Petrich} 
 {Inverse semigroups}
 {John Wiley \& Sons, New York, 1984.}

\i{[Pr1]}{G. Preston} 
 {Inverse semi-groups}
 {J. London Math. Soc. {\bf 29}, (1954). 396--403.}

\i{[Pr2]}{G. Preston} 
 {Monogenic inverse semigroups}
 {J. Austral. Math. Soc. Ser. A {\bf 40} (1986), no. 3, 321--342.}

\i{[SZ]}{Z. Shan, P.T. Zhu}
 {On $(a,b,k)$-partitions of positive integers}
 {Southeast Asian Bull. Math. {\bf 17} (1993), no. 1, 51--58.}

\i{[Vag]}{V. Vagner}
 {Generalized groups}
 {Doklady Akad. Nauk SSSR (N.S.) {\bf 84}, (1952). 1119--1122.}

\i{[ZC]}{P. Zhu, M. Chen}
 {$(a,b,k)$ $(k<0)$-additive partition of positive integers}
 {Chinese Sci. Bull. {\bf 41} (1996), no. 24, 2029--2037.}

\uj
\noindent
{\it Email:} {\tt nandor.sieben@nau.edu}

\noindent
{\nagy Department of Mathematics, Northern Arizona 
University, Flagstaff, AZ 86011-5717}

\end